\documentclass[5p,numbers,sort&compress,intlimits]{elsarticle}
\usepackage[T1]{fontenc}
\usepackage[utf8]{inputenc}
\usepackage{color}
\usepackage{mathtools}
\usepackage{amsmath}
\usepackage{amsthm}
\usepackage{amssymb}

\makeatletter
\theoremstyle{plain}
\newtheorem{thm}{\protect\theoremname}
\theoremstyle{plain}
\newtheorem{lem}[thm]{\protect\lemmaname}

\makeatother

\providecommand{\lemmaname}{Lemma}
\providecommand{\theoremname}{Theorem}

\begin{document}

\begin{frontmatter}{}

\title{Some results on thermopiezoelectricity of nonsimple materials}

\author[1]{Michele Ciarletta}

\ead{mciarletta@unisa.it}

\author[2]{Martina Nunziata}

\ead{mnunziata@unisa.it}

\author[2]{Francesca Passarella}

\ead{fpassarella@unisa.it}

\author[2]{Vincenzo Tibullo\corref{cor1}}

\ead{vtibullo@unisa.it}

\cortext[cor1]{Corresponding author}

\affiliation[1]{organization={Dipartimento di Ingegneria Civile, Università di Salerno}}

\affiliation[2]{organization={Dipartimento di Matematica, Università di Salerno}}
\begin{abstract}
In this paper, we consider the linear theory for a model of a thermopiezoelectric
nonsimple material adopting the entropy production inequality proposed
by Green and Laws as presentated in \citep{Passarella2022}. We establish
reciprocity theorems and a variational principle for homogeneous and
anisotropic thermopiezoelectric nonsimple materials with a center
of symmetry. The proof of these theorems use the time convolution
product and an alternative formulation of the field equations. Moreover,
a uniqueness result is established without using the definiteness
assumptions on internal energy. 
\end{abstract}

\end{frontmatter}{}

\emph{This paper is dedicated to prof. Brian Straughan, a great researcher
but above all a good friend.}

\section{Introduction}

In \citet{Passarella2022}, the authors derived a theory for a thermopiezoelectric
body in which the second gradient of displacement field and the second
gradient of the electric potential are included in the set of independent
constitutive variables. They obtained the thermodynamic restrictions
and constitutive equations by using the entropy production inequality
proposed by \citet{Green1972a}.

The theory proposed by Green and Laws is one of the theories that
predict a finite velocity for the propagation of thermal signals (see
the reviews of \citet{Chandrasekharaiah1998,Chandrasekharaiah1986}).
As shown by~\citet{Iesan2004c}, they make use of an entropy inequality
in which a new constitutive function appears with the role of thermodynamic
temperature (see e.g.~\citet{Passarella2013}). In addition to the
finite velocity of heat waves, this theory also results in a symmetric
heat conductivity tensor.

By the other end, the origin of the theory of nonsimple elastic materials
goes back to the works of \citet{Toupin1962,Toupin1964}, and \citet{Mindlin1964}.
\citet{Toupin1965} applied the general theory of materials of grade
2 to the problem of surface deformations of a crystal. They showed
that initial stress and hyperstress in a uniform crystal give rise
to a deformation of a thin boundary layer near a free surface such
as that observed in electron diffraction experiments. Strain gradient
theory of thermoelasticity was first presented in~\citet{Ahmadi1975,Batra1976}.
The gradient theory of elasticity becomes important because it is
adequate to investigate problems related to size effects and nanotechnology.
In the regime of micron and nano-scales, experimental evidence and
observations have suggested that classical continuum theories do not
suffice for an accurate and detailed description of corresponding
deformation phenomena. The theory of nonsimple thermoelastic materials
has been discussed in various papers (see for example~\citep{Mindlin1968,Ahmadi1977,Iesan1983,Ciarletta1989a,Kalpakides2002,Aouadi2019,Aouadi2020,Aouadi2020a,Bartilomo1997}).

Furthermore, \citet{Kalpakides2002} have established a theory of
electroelasticity including both strain gradient and electric field
gradient. They report that taking into account of the second spatial
gradient of the motion makes sense especially in crack problems, moreover
taking into account of second gradient of the electrical potential
implies the presence of quadrupole polarization into the continuum
model, of practical interest for problems concerning surface effects.

The problem of the interaction of the electromagnetic field with the
motion of elastic solids was the subject of important investigations
(see e.g.~\citep{Eringen2004,Truesdell1960,Parkus1972,Grot1976,Nowacki1983,Maugin1988,Morro1991}
and the literature cited therein). Certain crystals (for example quartz)
when subject to stress, become electrically polarized (piezoelectric
effect). Conversely, an external electromagnetic field can produce
deformation in a piezoelectric crystal.

In section 2, we begin by summarizing the fundamental equations based
on the linear theory of thermopiezoelectric nonsimple materials as
established in \citet{Passarella2022} and we consider in particular
the case of center-symmetric materials. In section 3, we define a
mixed initial-boundary value problem under non-homogeneous initial
conditions and present a characterization of the mixed initial-boundary
value problem in an alternative way, by including the initial conditions
into the field equations. In section 4, starting from a reciprocity
relation which involves two processes at different times, a uniqueness
result is established without using the definiteness assumptions on
internal energy. Moreover, a reciprocity theorem is presented. In
sections 5 and 6, another reciprocity theorem based on the convolution
product and a variational principle are derived (see also \citep{Passarella2011a}).

\section{Basic equations}

We consider a body that at some instant occupies the region $B$ of
the Euclidean three-dimensional space and is bounded by the smooth
surface $\partial B$. The motion of the body is referred to the reference
configuration $B$ and to a fixed system of rectangular Cartesian
axes $Ox_{i}$ $(i=1,2,3)$.

We shall employ the usual summation and differentiation conventions:
Latin subscripts are understood to range over the integers $(1,2,3)$,
summation over repeated subscripts is implied and subscripts preceded
by a comma denote partial differentiation with respect to the corresponding
Cartesian coordinate. In what follows we use a superposed dot to denote
partial differentiation with respect to the time $t$.

As already done by \citet{Passarella2022}, we consider the linear
theory for a model of a thermopiezoelectric nonsimple materials adopting
the entropy production inequality proposed by \citet{Green1972a},
in local form 
\[
\rho\phi\dot{\eta}\geq\rho h-q_{i,i}+\frac{1}{\phi}q_{i}\phi_{,i}
\]
where $\rho>0$ is the reference mass density, $\eta$ is the entropy
per unit mass, $q_{i}$ is the heat flux vector, $h$ is the heat
supply per unit mass and unit time. The function $\phi$ is a new
strictly positive thermal function.

Let $u_{i}$ be the displacement vector, $\theta$ the difference
of the absolute temperature $T$ from the absolute temperature in
the reference configuration $T_{0}>0$ (i.e. $\theta=T-T_{0}$), $D_{i}$
the electric displacement vector field, $E_{i}$ the electric vector
field, then the linear theory is governed by the following local balance
equations defined in $B\times(0,\infty)$: 
\begin{align}
 & \tau_{ji,j}-\mu_{kji,kj}+\rho f_{i}=\rho\ddot{u}_{i},\label{eq:motion}\\
 & \sigma_{i,i}-Q_{ji,ji}=g,\label{eq:electric-field}\\
 & \rho T_{0}\dot{\eta}+q_{j,j}-\rho h=0,\label{eq:energy}
\end{align}
with $\tau_{ji}=t_{ji}+\mu_{kji,k},$ $\sigma_{i}=D_{i}+Q_{ji,j}$
and 
\begin{alignat}{2}
 & e_{ij}=\frac{1}{2}(u_{i,j}+u_{j,i}),\; &  & \kappa_{ijk}=u_{k,ij},\;\label{eq:geometric-1}\\
 & \beta_{i}=\theta_{,i},\ E_{i}=-\varphi_{,i},\; & \  & V_{ij}=-\varphi_{,ij}.\label{eq:geometric-2}
\end{alignat}
Here, $t_{ij}$ is the stress tensor, $\mu_{kji}$ the hyperstress
tensor, $f_{i}$ the external body force per unit mass, $Q_{ji}$
the electric quadrupole, $g$ the density of free charges and $\varphi$
the electric potential. 

We restrict our attention to a homogeneous, center-symmetric material,
so that the constitutive equations defined on $\bar{B}\times I$ (with
$I=[0,\infty)$) are 
\begin{equation}
\begin{aligned}\tau_{ij} & =a_{ijkl}^{(11)}e_{kl}+a_{ijkl}^{(17)}V_{kl}+a_{ij}^{(14)}\left(\theta+\beta\dot{\theta}\right),\\
\mu_{ijk} & =a_{ijklhm}^{(22)}\kappa_{lhm}+a_{ijkl}^{(23)}E_{l},\\
-\sigma_{i} & =a_{jkli}^{(23)}\kappa_{jkl}+a_{ij}^{(33)}E_{j},\\
-Q_{ij} & =a_{klij}^{(17)}e_{kl}+a_{ijkl}^{(77)}V_{kl}+a_{ij}^{(47)}\left(\theta+\beta\dot{\theta}\right),\\
-\rho\eta & =a_{ij}^{(14)}e_{ij}+a_{ij}^{(47)}V_{ij}+c\left(\theta+\beta'\dot{\theta}\right),\\
-\frac{q_{i}}{T_{0}} & =k_{ij}\beta_{j},
\end{aligned}
\label{eq:constitutive}
\end{equation}
with $\beta'=\beta+\gamma/(c\beta)$, $\beta\neq0$, $c\neq0$. We
assume that the coefficients in~\eqref{eq:constitutive} satisfy
the following symmetry relations
\begin{equation}
\begin{alignedat}{2} & a_{ijkl}^{(11)}=a_{jikl}^{(11)}=a_{klij}^{(11)}, & \quad & a_{ij}^{(14)}=a_{ji}^{(14)},\\
 & a_{ijklhm}^{(22)}=a_{jiklhm}^{(22)}=a_{lhmijk}^{(22)}, &  & a_{ijkl}^{(23)}=a_{jikl}^{(23)},\\
 & a_{ij}^{(33)}=a_{ji}^{(33)},\ \quad k_{ij}=k_{ji} &  & a_{ijkl}^{(17)}=a_{jikl}^{(17)},\\
 & a_{ijkl}^{(77)}=a_{jikl}^{(77)}=a_{klij}^{(77)}, &  & a_{ij}^{(47)}=a_{ji}^{(47)}.
\end{alignedat}
\label{eq:symmetries}
\end{equation}
Furthermore, it is 
\[
\begin{alignedat}{3} & e_{ij}=e_{ji}, & \  & \kappa_{ijk}=\kappa_{jik}, & \  & V_{ij}=V_{ji},\\
 & \tau_{ji}=\tau_{ij}, &  & \mu_{kji}=\mu_{jki}, &  & Q_{ij}=Q_{ji}.
\end{alignedat}
\]
The dissipation inequality implies that the quadratic form $\mathcal{P}$
is positive semi-definite, i.e. 
\begin{equation}
\mathcal{P}(\xi,\eta_{i})=\frac{\gamma}{\beta}\xi^{2}+k_{ij}\eta_{i}\eta_{j}\geq0,\qquad\forall\xi,\eta_{i}.\label{eq:dissipation-quadratic-form}
\end{equation}
The inequality~\eqref{eq:dissipation-quadratic-form} is equivalent
to 
\begin{equation}
\frac{\gamma}{\beta}\geq0,\quad k_{ij}\eta_{i}\eta_{j}\geq0,\qquad\forall\eta_{i}.\label{eq:semi-positivity-conditions}
\end{equation}
It results that the tensor $k_{ij}$ is positive semi-definite.

\section{Mixed initial-boundary value problem }

Now, we denote with $\Pi$ the mixed initial-boundary value problem
defined by eqs.~\eqref{eq:motion}-\eqref{eq:constitutive} with
the restriction~\eqref{eq:semi-positivity-conditions}, the following
initial conditions 
\begin{equation}
u_{i}(0)=u_{i}^{0},\quad\dot{u}(0)=v_{i}^{0},\quad\theta(0)=\theta^{0},\quad\eta(0)=\eta^{0},\label{eq:initial-conditions}
\end{equation}
in $\bar{B}$ and the following boundary conditions 
\begin{equation}
\begin{alignedat}{5} & u_{i} &  & =\hat{u}_{i} & \ \text{on}\;S_{1}\times I,\quad & P_{i} &  & =\hat{P}_{i} & \ \text{on}\;\Sigma_{1}\times I,\\
 & \mathcal{D}u_{i} &  & =\hat{d}_{i} & \text{on}\;S_{2}\times I,\quad & R_{i} &  & =\hat{R}_{i} & \text{on}\;\Sigma_{2}\times I,\\
 & \theta &  & =\hat{\theta} & \text{on}\;S_{3}\times I,\quad & q &  & =\hat{q} & \text{on}\;\Sigma_{3}\times I,\\
 & \varphi &  & =\hat{\varphi} & \text{on}\;S_{4}\times I,\quad & \Lambda &  & =\hat{\Lambda} & \text{on}\;\Sigma_{4}\times I,\\
 & \mathcal{D}\varphi &  & =\hat{\omega} & \text{on}\;S_{5}\times I,\quad & H &  & =\hat{H} & \text{on}\;\Sigma_{5}\times I,
\end{alignedat}
\label{eq:boundary-conditions}
\end{equation}
where, denoted by $n_{i}$ the outward unit normal vector to the boundary
surface $\partial B$, $q$ is the heat flux, i.e. $q=q_{i}n_{i}$,
and $\{S_{i},\Sigma_{i}\}$ are a subset of $\partial B$ such that,
considering the closure relative to $\partial B$, 
\[
\bar{S}_{i}\cup\Sigma_{i}=\partial B\quad S_{i}\cap\Sigma_{i}=\varnothing,\qquad i=1,\ldots,5,
\]
and we have \citep{Iesan2004c} {\small{}
\[
\begin{aligned} & P_{i}=(\tau_{ji}-\mu_{kji,k})n_{j}-\mathcal{D}_{j}(\mu_{kji}n_{k})+(\mathcal{D}_{l}n_{l})\mu_{kji}n_{k}n_{j}\\
 & \Lambda=(\sigma_{j}-Q_{kj,k})n_{j}-\mathcal{D}_{j}(Q_{kj}n_{k})+(\mathcal{D}_{l}n_{l})Q_{kj}n_{k}n_{j},
\end{aligned}
\]
}and 
\[
R_{i}=\mu_{kji}n_{k}n_{j},\qquad H=Q_{kj}n_{k}n_{j},
\]
where $\mathcal{D}\equiv n_{i}\partial/\partial x_{i}$ is the normal
derivative operator and $\mathcal{D}_{i}\equiv(\delta_{ij}-n_{i}n_{j})\partial/\partial x_{j}$
the surface gradient operator. We can prove that the functions $P_{i}$,
$R_{i}$, $\Lambda$ and $H$ are such that, for different times $r,s\in I$
\begin{equation}
\begin{aligned}\int_{\partial B} & \left[\left(\tau_{ji}(r)-\mu_{kji,k}(r)\right)u_{i}(s)+\mu_{jli}(r)u_{i,l}(s)\right]n_{j}da\\
 & \ =\int_{\partial B}\left[P_{i}(r)u_{i}(s)+R_{i}(r)\mathcal{D}u_{i}(s)\right]da,\\
\int_{\partial B} & \left[\left(\sigma_{j}(r)-Q_{ij,i}(r)\right)\varphi(s)+Q_{ji}(r)\varphi_{,i}(s)\right]n_{j}da\\
 & =\int_{\partial B}\left[\Lambda(r)\varphi(s)+H(r)\mathcal{D}\varphi(s)\right]da.
\end{aligned}
\label{eq:surface-integral-identities}
\end{equation}
All right-hand terms in eqs. \eqref{eq:initial-conditions} and \eqref{eq:boundary-conditions},
along with $f_{i}$, $g$ and $h$ are the given data of the considered
mixed initial-boundary value problem $\Pi$ and are prescribed continuous
functions. We denote the given data by 
\[
\Gamma=\left(f_{i},g,h,u_{i}^{0},v_{i}^{0},\theta^{0},\eta^{0},\hat{u}_{i},\hat{d}_{i},\hat{\theta},\hat{\varphi},\hat{\omega},\hat{P}_{i},\hat{R}_{i},\hat{q},\hat{\Lambda},\hat{H}\right).
\]
Let us define an ordered array of functions 
\[
\pi=(u_{i},\theta,\varphi,e_{ij},\kappa_{ijk},\beta_{i},E_{i},V_{ij},\tau_{ij},\mu_{ijk},\sigma_{i},Q_{ji},\eta,q_{i})
\]
as an admissible process on $\bar{B}\times I$ with the following
properties 
\begin{enumerate}
\item \begin{flushleft}
$u_{i}\in C^{4,2}(\bar{B}\times I)$, $\varphi\in C^{4,0}(\bar{B}\times I)$,
$\theta\in C^{2,2}(\bar{B}\times I)$, \\
$e_{ij},V_{ij}\in C^{2,1}(\bar{B}\times I)$, $\eta\in C^{0,1}(\bar{B}\times I),$\\
$\kappa_{ijk},E_{i},\mu_{ijk},Q_{ij}\in C^{2,0}(\bar{B}\times I)$,\\
$\beta_{i},\tau_{ij},q_{i}\in C^{1,0}(\bar{B}\times I)$; 
\par\end{flushleft}
\item \begin{flushleft}
$\text{\ensuremath{e_{ij}}=\ensuremath{e_{ji}}}$, $\kappa_{ijk}=\kappa_{jik}$,
$\,V_{ij}=V_{ji}$, $\tau_{ji}=\tau_{ij}$, \\
$\mu_{kji}=\mu_{jki}$, $Q_{ji}=Q_{ij}$ on $\bar{B}\times I$. 
\par\end{flushleft}

\end{enumerate}
We say that $\pi$ is a process corresponding to the supply terms
$(f_{i},g,h)$ if $\pi$ is an admissible process that satisfies the
fundamental system of field equations~\eqref{eq:motion}-\eqref{eq:constitutive}
with the restriction~\eqref{eq:semi-positivity-conditions} on $B\times[0,\infty)$.

Then, if a process $\pi$ satisfies the initial conditions~\eqref{eq:initial-conditions}
and the boundary conditions~\eqref{eq:boundary-conditions}, we identify
it as a solution of the mixed initial-boundary value problem $\Pi$.

The set $\mathcal{\ensuremath{V}}$ of all admissible processes on
$\bar{B}\times I$ can be considered as a vector space. We denote
by $\mathcal{K}\subseteq\mathcal{\ensuremath{V}}$ the set of all
solution of the mixed initial-boundary value problem in concern.

Following \citet{Gurtin1964}, we will give an alternative formulation
of the problem \eqref{eq:motion}-\eqref{eq:constitutive} in which
the initial conditions~\eqref{eq:initial-conditions} are incorporated
into the field equations. To this aim, we introduce the product of
convolution as follows 
\[
\left[f_{1}*f_{2}\right](t)=\int_{0}^{t}f_{1}(\tau)f_{2}(t-\tau)d\tau\qquad\quad\forall t\in I
\]
for any two continuous functions $f_{1}$, $f_{2}$, on $\bar{B}\times I$.
It is useful to introduce
\begin{align*}
 & l(t)=1, &  & \xi(t)=\left[l*l\right](t)=t,\\
 & \zeta(t)=\frac{1}{\beta}e^{-t/\beta}, &  & \chi(t)=\frac{1}{\beta'}e^{-t/\beta'},
\end{align*}
with $\beta'\neq0$, and we note in particular that
\[
\left[l*f\right](t)=\int_{0}^{t}f(\tau)d\tau.
\]
From these definitions, we easily obtain 
\begin{equation}
\begin{aligned} & \xi*\ddot{u}_{i}=u_{i}-v_{i}^{0}t-u_{i}^{0},\\
 & l*\dot{\theta}=\theta-\theta^{0}\qquad l*\dot{\eta}=\eta-\eta^{0},\\
 & \zeta*(\theta+\beta\dot{\theta})=\theta-\theta^{0}e^{-t/\beta},\\
 & \chi*(\theta+\beta'\dot{\theta})=\theta-\theta^{0}e^{-t/\beta'},
\end{aligned}
\label{eq:convolution-properties}
\end{equation}
so that we can prove the following Lemma.
\begin{lem}
Let $\pi\in\mathcal{V}$, then, $\pi$ satisfies eqs.~\eqref{eq:motion},
\eqref{eq:electric-field}, \eqref{eq:energy}, \eqref{eq:constitutive}
and the initial conditions~\eqref{eq:initial-conditions} if and
only if 
\begin{equation}
\begin{aligned} & \xi*\left(\tau_{ji,j}-\mu_{kji,kj}\right)-\rho u_{i}+\mathcal{F}_{i}=0,\\
 & \xi*\left(\sigma_{i,i}-Q_{ji,ji}-g\right)=0,\\
 & \rho T_{0}\eta+l*q_{i,i}-\mathcal{H}=0,\\
 & +\zeta*\left(\tau_{ij}-\hat{\tau}_{ij}\right)-a_{ij}^{(14)}\theta+\mathcal{L}_{ij}=0,\\
 & +\zeta*\left(\mu_{ijk}-\hat{\mu}_{ijk}\right)=0,\\
 & -\zeta*\left(\sigma_{i}-\hat{\sigma}_{i}\right)=0,\\
 & -\zeta*\left(Q_{ij}-\hat{Q}_{ij}\right)-a_{ij}^{(47)}\theta+\mathcal{M}_{ij}=0,\\
 & -\chi*\left(\rho\eta-\rho\hat{\eta}\right)-c\theta+\mathcal{R}=0,\\
 & -\frac{1}{T_{0}}l*\left(q_{i}-\hat{q}_{i}\right)=0,
\end{aligned}
\label{eq:convolution-alternative}
\end{equation}
 where 
\begin{equation}
\begin{aligned}\hat{\tau}_{ij} & =a_{ijkl}^{(11)}e_{kl}+a_{ijkl}^{(17)}V_{kl},\\
\hat{\mu}_{ijk} & =a_{ijklhm}^{(22)}\kappa_{lhm}+a_{ijkl}^{(23)}E_{l},\\
-\hat{\sigma}_{i} & =a_{jkli}^{(23)}\kappa_{jkl}+a_{ij}^{(33)}E_{j},\\
-\hat{Q}_{ij} & =a_{klij}^{(17)}e_{kl}+a_{ijkl}^{(77)}V_{kl},\\
-\rho\hat{\eta} & =a_{ij}^{(14)}e_{ij}+a_{ij}^{(47)}V_{ij},\\
-\hat{q}_{i} & =T_{0}k_{ij}\beta_{j},
\end{aligned}
\label{eq:constitutive-reduced}
\end{equation}
and 
\begin{equation}
\begin{aligned} & \mathcal{F}_{i}=\rho\xi*f_{i}+\rho v_{i}^{0}t+\rho u_{i}^{0}, &  & \mathcal{H}=\rho T_{0}\eta^{0}+l*\rho h,\\
 & \mathcal{L}_{ij}=a_{ij}^{(14)}\theta^{0}e^{-t/\beta}, &  & \mathcal{M}_{ij}=a_{ij}^{(47)}\theta^{0}e^{-t/\beta},\\
 & \mathcal{R}=c\theta^{0}e^{-t/\beta'}.
\end{aligned}
\label{eq:convoluted-data-defs}
\end{equation}
\end{lem}

\begin{thm}
A process $\pi$ is a solution of the mixed initial-boundary value
problem in concern $\Pi$ if and only if it satisfies eqs.~\eqref{eq:geometric-1},
\eqref{eq:geometric-2}, \eqref{eq:convolution-alternative} with
the restriction~\eqref{eq:semi-positivity-conditions} and the boundary
conditions~\eqref{eq:boundary-conditions}. 
\end{thm}

\section{Uniqueness and reciprocity theorems}

We consider the body $B$ subjected to two different sets of external
data 
\[
\begin{aligned}\Gamma^{(\alpha)} & =\biggl(f_{i}^{(\alpha)},g^{(\alpha)},h^{(\alpha)},u_{i}^{0(\alpha)},v_{i}^{0(\alpha)},\theta^{0(\alpha)},\eta^{0(\alpha)},\hat{u}_{i}^{(\alpha)},\\
 & \hat{d}_{i}^{(\alpha)},\hat{\theta}^{(\alpha)},\hat{\varphi}^{(\alpha)},\hat{\omega}^{(\alpha)},\hat{P}_{i}^{(\alpha)},\hat{R}_{i}^{(\alpha)},\hat{\Lambda}^{(\alpha)},\hat{H}^{(\alpha)},\hat{q}^{(\alpha)}\biggr),
\end{aligned}
\]
with $\alpha=1,2$, and denote the corresponding solutions of the
mixed initial-boundary problem as 
\begin{align*}
\pi{}^{(\alpha)}=( & u_{i}^{(\alpha)},\theta{}^{(\alpha)},\varphi{}^{(\alpha)},e_{ij}^{(\alpha)},\kappa_{ijk}^{(\alpha)},\beta_{i}^{(\alpha)},E_{i}^{(\alpha)},V_{ij}^{(\alpha)},\\
 & \tau_{ij}^{(\alpha)},\mu_{ijk}^{(\alpha)},\sigma_{i}^{(\alpha)},Q_{ji}^{(\alpha)},\eta^{(\alpha)},q_{i}^{(\alpha)})\in\mathcal{K}.
\end{align*}
Moreover, corresponding to $\Gamma^{(\alpha)}$ and $\pi^{(\alpha)}$,
let's define $\mathcal{F}_{i}^{(\alpha)}$, $\mathcal{H}^{(\alpha)}$,
$\mathcal{L}_{ij}^{(\alpha)}$, $\mathcal{R}^{(\alpha)}$, $\mathcal{M}_{ij}^{(\alpha)}$,
$\hat{\tau}{}_{ij}^{(\alpha)}$, $\hat{\mu}{}_{ijk}^{(\alpha)}$,
$\hat{\sigma}_{i}^{(\alpha)}$, $\hat{Q}{}_{ij}^{(\alpha)}$, $\hat{q}_{i}^{(\alpha)}$
through eqs. \eqref{eq:constitutive-reduced} and \eqref{eq:convoluted-data-defs}.

Henceforth, the dependence on time will be explicit, while the\textcolor{red}{{}
}dependence on $\mathbf{x}$ will remain implicit. It is useful for
the following to introduce the functions 
\[
\begin{aligned}\hat{\mathcal{S}}_{\alpha\beta}(r,s)=\  & \hat{\tau}_{ji}^{(\alpha)}(r)e_{ij}^{(\beta)}(s)+\hat{\mu}_{kji}^{(\alpha)}(r)\kappa_{kji}^{(\beta)}(s)\\
- & \hat{\sigma}_{i}^{(\alpha)}(r)E_{i}^{(\beta)}(s)-\hat{Q}_{ji}^{(\alpha)}(r)V_{ij}^{(\beta)}(s)
\end{aligned}
\]
and 
\begin{equation}
\begin{aligned}\mathcal{S}_{\alpha\beta}(r,s)=\  & \tau_{ji}^{(\alpha)}(r)e_{ij}^{(\beta)}(s)+\mu_{kji}^{(\alpha)}(r)\kappa_{kji}^{(\beta)}(s)\\
- & \sigma_{i}^{(\alpha)}(r)E_{i}^{(\beta)}(s)-Q_{ji}^{(\alpha)}(r)V_{ij}^{(\beta)}(s)
\end{aligned}
\label{eq:S}
\end{equation}
for different times $r,s\in I$. We can prove 
\begin{equation}
\begin{aligned}\hat{\mathcal{S}}_{\alpha\beta}(r,s) & =\mathcal{\hat{S}}_{\beta\alpha}(s,r),\\
\mathcal{S}_{\alpha\beta}(r,s) & =\hat{\mathcal{S}}_{\alpha\beta}(r,s)-A\theta^{(\alpha)}(r)\rho\hat{\eta}^{(\beta)}(s),
\end{aligned}
\label{eq:S-relations}
\end{equation}
where $A$ is the following differential operator 
\begin{equation}
A=I+\beta\dfrac{\partial}{\partial t},\label{eq:A-operator}
\end{equation}
with $I$ the identity operator. If we define 
\begin{align*}
\mathcal{\hat{T}}_{\alpha\beta}(t) & =\xi*\int_{0}^{t}\hat{\mathcal{S}}_{\alpha\beta}(\tau,t-\tau)d\tau,\\
\mathcal{T}_{\alpha\beta}(t) & =\xi*\int_{0}^{t}\mathcal{S}_{\alpha\beta}(\tau,t-\tau)d\tau,
\end{align*}
eqs.~\eqref{eq:S-relations} imply 
\begin{equation}
\begin{alignedat}{1}\mathcal{\hat{T}}_{\alpha\beta}(t) & =\mathcal{\hat{T}}_{\beta\alpha}(t),\\
\mathcal{T}_{\alpha\beta}(t) & =\mathcal{\hat{T}}_{\alpha\beta}(t)-\xi*A\theta^{(\alpha)}*\rho\hat{\eta}^{(\beta)}(t).
\end{alignedat}
\label{eq:T-relations}
\end{equation}
For what follows it is useful to remark that, by eqs. \eqref{eq:convolution-properties}$_{3}$,
\eqref{eq:convoluted-data-defs}, \eqref{eq:S}, \eqref{eq:T-relations}$_{2}$,
the following relations hold 
\[
\begin{aligned}\mathcal{T}_{\alpha\beta}(t)=\xi*\biggl[ & \tau_{ji}^{(\alpha)}*e_{ij}^{(\beta)}(t)+\mu_{kji}^{(\alpha)}*\kappa_{kji}^{(\beta)}(t)\\
- & \sigma_{i}^{(\alpha)}*E_{i}^{(\beta)}(t)-Q_{ji}^{(\alpha)}*V_{ij}^{(\beta)}(t)\biggr],
\end{aligned}
\]
\begin{equation}
\begin{aligned}\zeta*\mathcal{T}_{\alpha\beta}(t)=\  & \mathcal{\zeta*\hat{T}}_{\alpha\beta}(t)-\xi*\theta^{(\alpha)}*\rho\hat{\eta}^{(\beta)}(t)\\
-\  & \xi*\mathcal{L}_{ij}^{(\alpha)}*e_{ij}^{(\beta)}(t)-\xi*\mathcal{M}_{ij}^{(\alpha)}*V_{ij}^{(\beta)}(t).
\end{aligned}
\label{eq:zeta-T}
\end{equation}
On the other hand, by eqs. \eqref{eq:geometric-1} and \eqref{eq:geometric-2}
we have 
\begin{equation}
\begin{aligned} & \mathcal{S}_{\alpha\beta}(r,s)=-\left[\tau_{kj,j}^{(\alpha)}(r)-\mu_{jik,ij}^{(\alpha)}(r)\right]u_{k}^{(\beta)}(s)\\
 & \;-\left[\sigma_{i,i}^{(\alpha)}(r)-Q_{ij,ij}^{(\alpha)}(r)\right]\varphi^{(\beta)}(s)\\
 & \;+\biggl[\left[\tau_{kj}^{(\alpha)}(r)-\mu_{jik,i}^{(\alpha)}(r)\right]u_{k}^{(\beta)}(s)+\mu_{ijk}^{(\alpha)}(r)u_{k,i}^{(\beta)}(s)\\
 & \;+\left[\sigma_{j}^{(\alpha)}(r)-Q_{ij,i}^{(\alpha)}(r)\right]\varphi^{(\beta)}(s)+Q_{ji}^{(\alpha)}(r)\varphi_{,i}^{(\beta)}(s)\biggr]_{,j},
\end{aligned}
\label{eq:S-surface}
\end{equation}
so that, taking into account that eqs. \eqref{eq:convolution-alternative}
hold for $\pi^{(1,2)}\in\mathcal{K}$, using eqs.\eqref{eq:surface-integral-identities}
and the divergence theorem, we arrive to 
\begin{equation}
\begin{aligned} & \int_{B}\mathcal{T}_{\alpha\beta}(t)dv=-\int_{B}\rho u_{i}^{(\alpha)}*u_{i}^{(\beta)}(t)dv\\
 & \ +\int_{B}\left[\mathcal{F}_{i}^{(\alpha)}*u_{i}^{(\beta)}(t)-\xi*g^{(\alpha)}*\varphi^{(\beta)}(t)\right]dv\\
 & \;+\int_{\partial B}\xi*\Bigl[P_{i}^{(\alpha)}*u_{i}^{(\beta)}(t)+R_{i}^{(\alpha)}*\mathcal{D}u_{i}^{(\beta)}(t)\\
 & \;+\Lambda^{(\alpha)}*\varphi^{(\beta)}(t)+H^{(\alpha)}*\mathcal{D}\varphi^{(\beta)}(t)\Bigr]da.
\end{aligned}
\label{eq:T-integral}
\end{equation}
In this section, we set
\[
\Theta_{,j}^{(\alpha)}(t)=\left[l*\theta_{,j}^{(\alpha)}\right](t)=\int_{0}^{t}\theta_{,j}^{(\alpha)}(\tau)d\tau.
\]
Now, we obtain the following reciprocity relation which involves two
processes at different times 
\begin{lem}
\label{thm:gamma-symmetry} Let $\pi^{(1,2)}\in\mathcal{K}$. Then
\begin{equation}
\Gamma_{\alpha\beta}(r,s)=\Gamma_{\beta\alpha}(s,r),\qquad\forall r,s\in I,\ \forall\alpha,\beta=1,2,\label{eq:gamma-symmetry}
\end{equation}
where we define 
\begin{equation}
\begin{aligned} & \Gamma_{\alpha\beta}(r,s)=\int_{B}\biggl[\rho f_{i}^{(\alpha)}(r)u_{i}^{(\beta)}(s)-g^{(\alpha)}(r)\varphi^{(\beta)}(s)\ \\
 & \ -\frac{1}{T_{0}}\mathcal{H}^{(\alpha)}(r)A\theta^{(\beta)}(s)\biggr]dv+\int_{\partial B}\biggl[P_{i}^{(\alpha)}(r)u_{i}^{(\beta)}(s)\\
 & \;+R_{i}^{(\alpha)}(r)\mathcal{D}u_{i}^{(\beta)}(s)+\Lambda^{(\alpha)}(r)\varphi^{(\beta)}(s)\\
 & \;+H^{(\alpha)}(r)\mathcal{D}\varphi^{(\beta)}(s)+\frac{1}{T_{0}}\left[l*q^{(\alpha)}\right](r)A\theta^{(\beta)}(s)\biggr]da\\
 & \ -\int_{B}\biggl[\rho\ddot{u}_{i}^{(\alpha)}(r)u_{i}^{(\beta)}(s)+\frac{1}{\beta}\gamma\dot{\theta}^{(\alpha)}(r)\theta^{(\beta)}(s)\\
 & \ -k_{ij}\Theta_{,j}^{(\alpha)}(r)A\theta_{,i}^{(\beta)}(s)\biggr]dv.
\end{aligned}
\label{eq:gamma}
\end{equation}
\end{lem}

\begin{proof}
The first step is to introduce the following function 
\begin{equation}
\begin{aligned} & J_{\alpha\beta}(r,s)=\mathcal{S}_{\alpha\beta}(r,s)-\rho\eta^{(\alpha)}(r)A\theta^{(\beta)}(s).\end{aligned}
\label{eq:J-step-1}
\end{equation}
Taking into account the constitutive equations \eqref{eq:constitutive}$_{3}$
and eqs. \eqref{eq:S-relations}, we have 
\begin{equation}
\begin{aligned} & J_{\alpha\beta}(r,s)-\frac{1}{\beta}\gamma\dot{\theta}^{(\alpha)}(r)\theta^{(\beta)}(s)=\\
 & =\hat{\mathcal{S}}_{\alpha\beta}(r,s)+cA\theta^{(\alpha)}(r)A\theta^{(\beta)}(s)\\
 & \;-\left(\rho\hat{\eta}^{(\alpha)}(r)A\theta^{(\beta)}(s)+A\theta^{(\alpha)}(r)\rho\hat{\eta}{}^{(\beta)}(s)\right)\;\\
 & \;+\gamma\dot{\theta}^{(\alpha)}(r)\dot{\theta}^{(\beta)}(s)\\
 & \;=J_{\beta s}(s,r)-\frac{1}{\beta}\gamma\dot{\theta}^{(\beta)}(s)\theta^{(\alpha)}(r).
\end{aligned}
\label{eq:J-step-2}
\end{equation}
On the other hand, eqs. \eqref{eq:convolution-alternative}$_{3}$,
\eqref{eq:S-surface} and \eqref{eq:J-step-1} lead to 
\begin{equation}
\mathcal{\int_{B}}\left[J_{\alpha\beta}(r,s)-\frac{1}{\beta}\gamma\dot{\theta}^{(\alpha)}(r)\theta^{(\beta)}(s)\right]dv=\Gamma_{\alpha\beta}(r,s),\label{eq:J-step-3}
\end{equation}
and consequently, we arrive to the desired result by \eqref{eq:J-step-2}
and \eqref{eq:J-step-3}. 
\end{proof}
We will use this Lemma to establish a uniqueness theorem with no definiteness
assumption on internal energy and a reciprocity theorem.

\subsection{Uniqueness theorem}

Let $\pi^{(1)}\in\mathcal{\ensuremath{K}}$ and call it $\pi$ for
simplicity, we take in eq. \eqref{eq:gamma} $r=t+\tau$ e $s=t-\tau$
with $\alpha=\beta=1$ and integrating from $0$ to $t$ we obtain
\begin{equation}
\begin{aligned} & \int_{0}^{t}\Gamma_{11}(t+\tau,t-\tau)d\tau=\int_{0}^{t}E(t+\tau,t-\tau)d\tau\\
 & -\int_{B}\int_{0}^{t}\biggl[\rho\ddot{u}_{i}(t+\tau)u_{i}(t-\tau)+\frac{1}{\beta}\gamma\dot{\theta}(t+\tau)A\theta(t-\tau)\\
 & \ -k_{ji}\Theta(t+\tau)A\theta_{,i}(t-\tau)\biggr]d\tau dv,
\end{aligned}
\label{eq:gamma-integral}
\end{equation}
where 
\[
\begin{aligned}E(r,s) & =\int_{B}\biggl[\rho f_{i}(r)u_{i}(s)-g(r)\varphi(s)-\frac{1}{T_{0}}\mathcal{H}(r)A\theta(s)\biggr]dv\ \\
 & \;+\int_{\partial B}\biggl[P_{i}(r)u_{i}(s)+R_{i}(r)\mathcal{D}u_{i}(s)+\Lambda(r)\varphi(s)\\
 & \;+H(r)\mathcal{D}\varphi(s)+\frac{1}{T_{0}}\left[l*q\right](r)A\theta(s)\biggr]da.
\end{aligned}
\]
Obviously, eq. \eqref{eq:gamma-symmetry} implies
\begin{equation}
\int_{0}^{t}\left[\Gamma_{11}(t+\tau,t-\tau)-\Gamma_{11}(t-\tau,t+\tau)\right]d\tau=0.\label{eq:gamma-difference}
\end{equation}
Let's use eqs. \eqref{eq:gamma-integral}, \eqref{eq:gamma-difference}
and the following relations 
\[
\begin{aligned} & \int_{0}^{t}\left[\ddot{u}_{i}(t+\tau)u_{i}(t-\tau)-\ddot{u}_{i}(t-\tau)u_{i}(t+\tau)\right]d\tau=\\
 & \qquad=\dot{u}_{i}(2t)u_{i}^{0}+u_{i}(2t)v_{i}^{0}-2u_{i}(t)\dot{u}_{i}(t),
\end{aligned}
\]
\[
\begin{aligned} & \int_{0}^{t}\left[\dot{\theta}(t+\tau)\theta(t-\tau)-\dot{\theta}(t-\tau)\theta(t+\tau)\right]d\tau=\\
 & \qquad=\theta(2t)\theta^{0}-\theta(t)\theta(t),
\end{aligned}
\]
{\small{}
\[
\begin{aligned}\int_{0}^{t} & \left[\Theta_{,i}(t-\tau)A\dot{\Theta}_{,j}(t+\tau)-\Theta_{,i}(t+\tau)A\dot{\Theta}_{,j}(t-\tau)\right]d\tau=\\
 & \negmedspace=-\Theta_{,i}(t)\Theta_{,j}(t)+\beta\left[\Theta_{,j}(2t)\theta_{,i}^{0}-2\Theta_{,j}(t)\theta_{,i}(t)\right],
\end{aligned}
\]
} to arrive to 
\begin{equation}
\begin{aligned}\int_{0}^{t} & \left[E(t+\tau,t-\tau)-E(t-\tau,t+\tau)\right]d\tau\\
 & \ -\int_{B}\rho\biggl\{\left[\dot{u}_{i}(2t)u_{i}^{0}+u_{i}(2t)v_{i}^{0}-2u_{i}(t)\dot{u}_{i}(t)\right]\\
 & \ +\frac{1}{\beta}\gamma\left[\theta(2t)\theta^{0}-\theta^{2}(t)\right]+k_{ij}\Bigl[-\Theta_{,i}(t)\Theta_{,j}(t)\\
 & \ +\beta\left[\Theta_{,j}(2t)\theta_{,i}^{0}-2\Theta_{,j}(t)\theta_{,i}(t)\right]\Bigr]\biggr\} dv=0.
\end{aligned}
\label{eq:E-difference}
\end{equation}
Eq.~\eqref{eq:E-difference} implies 
\begin{equation}
\begin{aligned} & \dot{G}(t)=-\int_{0}^{t}\left[E(t+\tau,t-\tau)-E(t-\tau,t+\tau)\right]d\tau\\
 & \ \int_{B}\biggl[\rho\left[\dot{u}_{i}(2t)u_{i}^{0}+u_{i}(2t)v_{i}^{0}\right]+\frac{1}{\beta}\gamma\theta(2t)\theta^{0}\\
 & +\beta k_{ij}\Theta_{,j}(2t)\theta_{,i}^{0}\biggr]dv,
\end{aligned}
\label{eq:G-derivative}
\end{equation}
where 
\begin{equation}
\begin{aligned} & G(t)=\int_{0}^{t}\int_{B}\mathcal{\mathcal{P}}\left[\theta(\tau),\Theta{}_{,i}(\tau)\right]dvd\tau\\
 & \quad+\int_{B}\left[\rho u_{i}(t)u_{i}(t)+\beta k_{ji}\Theta_{,j}(t)\Theta_{,i}(t)\right]dv,
\end{aligned}
\label{eq:G}
\end{equation}
with $\mathcal{\mathcal{P}}$ defined by eq.~\eqref{eq:dissipation-quadratic-form}. 

Now, we can prove the following uniqueness theorem 
\begin{thm}
Assume that 
\begin{enumerate}
\item $\beta$, $\gamma$ are strictly positive, 
\item the following quadratic form is definite
\[
F=a_{ji}^{(33)}E_{j}E_{i}+a_{jikl}^{(77)}V_{kl}V_{ji}.
\]
\end{enumerate}
If $S_{4}$ is nonempty, the initial-boundary values problem $\Pi$
has at most one solution. 
\end{thm}

\begin{proof}
Clearly, the difference $\pi$ of any two solutions of $\Pi$ corresponds
to null data. For this solution $\pi$, the function $G(t)$ defined
by eq.~\eqref{eq:G} vanishes initially and its derivative~\eqref{eq:G-derivative}
is identically zero, then $G(t)=0$ for all $t\in I.$ Since $\rho,\beta,\gamma>0$
and $\mathcal{P}$ is positive semi-definite, then for all $t\in I$
\begin{equation}
\begin{aligned} & u_{i}=0,\quad\theta=0\qquad\text{on}\;B\times I.\end{aligned}
\label{eq:u-theta-null}
\end{equation}
From eqs.~\eqref{eq:geometric-1}$_{1,2}$ it follows 
\begin{equation}
e_{ij}=0,\quad\kappa_{ijk}=0,\qquad\text{on}\;B\times I.\label{eq:e-kappa-null}
\end{equation}
Moreover, the constitutive equations~\eqref{eq:constitutive}$_{6}$
and the equation of energy~\eqref{eq:energy} with homogeneous initial
conditions imply 
\[
q_{i}=0,\quad\rho\eta=0\qquad\text{on}\;B\times I.
\]
On the other hand, it follows from eqs.~\eqref{eq:geometric-2},
\eqref{eq:constitutive}$_{3,4}$, \eqref{eq:u-theta-null}, \eqref{eq:e-kappa-null}
\begin{equation}
\begin{aligned}\int_{B} & F(t)dv=\!\int_{B}\Bigl[\sigma_{j}(t)\varphi_{,j}(t)+Q_{ij}(t)\varphi_{,ij}(t)\Bigr]dv.\end{aligned}
\label{eq:F-integral}
\end{equation}
Taking into account eqs.~\eqref{eq:electric-field}, \eqref{eq:surface-integral-identities}
and \eqref{eq:F-integral}, the divergence theorem and the null data,
we have 
\[
\begin{aligned}\int_{B}F(t)dv= & -\int_{B}g(t)\varphi(t)dv\\
 & +\int_{\partial B}\left[\Lambda(t)\varphi(t)+H(t)\mathcal{D}\varphi(t)\right]da=0.
\end{aligned}
\]
Consequently, given that $F$ is definite, we arrive to 

\begin{equation}
F(t)=0\quad\Rightarrow\quad E_{j}=0,\;V_{ji}=0,\quad\text{on}\;B\times I.\label{eq:E-V-null}
\end{equation}
Now, using eqs.~\eqref{eq:geometric-2}, \eqref{eq:constitutive}
and \eqref{eq:E-V-null}, we get 
\[
\begin{alignedat}{3} & \tau_{ij}=0, & \quad & \mu_{ijk}=0, & \quad\\
 & \sigma_{i}=0, & \quad & Q_{ji}=0, & \quad & \varphi=\text{const,}
\end{alignedat}
\]
on $B\times I$. If $S_{4}$ is nonempty, then 
\[
\varphi=0,\quad\text{on }S_{4}\times I\quad\Rightarrow\quad\varphi=0,\quad\text{on}\;B\times I
\]
and the proof is complete. 
\end{proof}

\subsection{Reciprocity theorem}

In this subsection we derive a reciprocity theorem based on Lemma
\ref{thm:gamma-symmetry} and following the method shown by~\citet{Iesan2004c} 
\begin{lem}
Let be $\pi^{(1,2)}\in\mathcal{K}$. Then we have 
\begin{equation}
I_{\alpha\beta}(t)=I_{\beta\alpha}(t),\qquad\forall t\in I,\ \forall\alpha,\beta=1,2,\label{eq:I-symmetry}
\end{equation}
with 
\begin{align*}
I_{\alpha\beta}(t)= & \int_{B}\biggl[\mathcal{F}_{i}^{(\alpha)}*u_{i}^{(\beta)}(t)-\xi*g^{(\alpha)}*\varphi^{(\beta)}(t)\biggr]dv\\
+ & \xi*\int_{\partial B}\biggl[P_{i}^{(\alpha)}*u_{i}^{(\beta)}(t)+R_{i}^{(\alpha)}*\mathcal{D}u_{i}^{(\beta)}(t)\\
+ & \Lambda^{(\alpha)}*\varphi^{(\beta)}(t)+H^{(\alpha)}*\mathcal{D}\varphi^{(\beta)}(t)\\
+ & \frac{1}{T_{0}}l*q^{(\alpha)}*A\theta^{(\beta)}(t)\biggr]da+\int_{B}\biggl[{\cal L}^{(\alpha)}*\theta^{(\beta)}(t)\\
- & \frac{1}{T_{0}}l*\mathcal{H}^{(\alpha)}*A\theta^{(\beta)}(t)+{\cal L}_{j}^{(\alpha)}*\theta_{,j}^{(\beta)}(t)\biggr]dv,
\end{align*}
where we define $\mathcal{F}_{i}^{(\alpha)}$, \textup{$\mathcal{H}^{(\alpha)}$,}
${\cal R}^{(\alpha)}$ by eq.~\eqref{eq:convoluted-data-defs} and
\begin{equation}
{\cal L}^{(\alpha)}=-\frac{1}{\beta}\gamma\theta^{0(\alpha)}\xi,\qquad{\cal L}_{j}^{(\alpha)}=\beta\xi*k_{ij}\theta_{,i}^{0(\alpha)}l.\label{eq:other-data-defs}
\end{equation}
\end{lem}

\begin{proof}
Taking into account eqs. \eqref{eq:convolution-properties}, \eqref{eq:convoluted-data-defs},
\eqref{eq:other-data-defs} and that 
\[
\begin{aligned} & k_{ij}\Theta_{,j}^{(\alpha)}*A\theta_{,i}^{(\beta)}(t)=k_{ij}\theta_{,j}^{(\alpha)}*l*\theta_{,i}^{(\beta)}(t)\\
 & \ +\beta k_{ij}\theta_{,j}^{(\alpha)}*\theta_{,i}^{(\beta)}(t)-\beta k_{ij}\theta_{,j}^{(\alpha)}*\theta_{,i}^{0(\beta)}(t),
\end{aligned}
\]
eq. \eqref{eq:gamma} leads to 
\[
\begin{aligned} & \xi*\int_{0}^{t}\Gamma_{\alpha\beta}(\tau,t-\tau)d\tau=\int_{B}\biggl[-\rho u_{i}^{(\alpha)}*u_{i}^{(\beta)}(t)\\
 & \ -\frac{1}{\beta}\gamma l*\theta^{(\alpha)}*\theta^{(\beta)}(t)+l*k_{ij}\Theta_{,j}^{(\alpha)}*\Theta_{,i}^{(\beta)}(t)\\
 & \;+\beta k_{ij}\Theta_{,j}^{(\alpha)}*\Theta_{,i}^{(\beta)}(t)+\mathcal{F}_{i}^{(\alpha)}*u_{i}^{(\beta)}(t)\\
 & \;-\xi*g^{(\alpha)}*\varphi^{(\beta)}(t)-{\cal L}^{(\alpha)}*\theta^{(\beta)}(t)\\
 & \ -\frac{1}{T_{0}}\xi*\mathcal{H}^{(\alpha)}*A\theta^{(\beta)}(t)-{\cal L}_{j}^{(\beta)}*\theta_{,j}^{(\alpha)}(t)\biggr]dv+\\
 & \;+\int_{\partial B}\xi*\biggl[P_{i}^{(\alpha)}*u_{i}^{(\beta)}(t)+R_{i}^{(\alpha)}*\mathcal{D}u_{i}^{(\beta)}(t)+\\
 & \;+\Lambda^{(\alpha)}*\varphi^{(\beta)}(t)+H^{(\alpha)}*\mathcal{D}\varphi^{(\beta)}(t)\\
 & \;+\frac{1}{T_{0}}l*q^{(\alpha)}*A\theta^{(\beta)}(t)\biggr]da.
\end{aligned}
\]
From this expression and Lemma \ref{thm:gamma-symmetry} it is easy
to prove that the following relation holds {\small{}
\begin{align*}
 & \xi*\int_{0}^{t}\Gamma_{\alpha\beta}(\tau,t-\tau)d\tau+\int_{B}\left[{\cal L}_{j}^{(\beta)}*\theta_{,j}^{(\alpha)}(t)+{\cal L}_{j}^{(\alpha)}*\theta_{,j}^{(\beta)}\right]dv\\
 & =\xi*\int_{0}^{t}\Gamma_{\beta\alpha}(\tau,t-\tau)d\tau+\int_{B}\left[{\cal L}_{j}^{(\beta)}*\theta_{,j}^{(\alpha)}(t)+{\cal L}_{j}^{(\alpha)}*\theta_{,j}^{(\beta)}\right]dv,
\end{align*}
}and this is equivalent to eq. \eqref{eq:I-symmetry}.
\end{proof}

\section{Alternative reciprocity theorem}

We now prove an alternative reciprocity theorem in which the operator
$A$ defined in eq. \eqref{eq:A-operator} is not used.
\begin{thm}
If we define {\small{}
\begin{equation}
\begin{aligned}\mathcal{I}_{\alpha\beta}(t)= & \chi*\zeta*\int_{B}\biggl[\mathcal{F}_{i}^{(\alpha)}*u_{i}^{(\beta)}(t)-\xi*g^{(\alpha)}*\varphi^{(\beta)}(t)\biggr]dv\\
+ & \chi*\xi*\int_{\partial B}\biggl\{\zeta*\biggl[P_{i}^{(\alpha)}*u_{i}^{(\beta)}(t)+R_{i}^{(\alpha)}*\mathcal{D}u_{i}^{(\beta)}(t)\\
+ & \Lambda^{(\alpha)}*\varphi^{(\beta)}(t)+H^{(\alpha)}*\mathcal{D}\varphi^{(\beta)}(t)\biggr]\\
+ & \frac{1}{T_{0}}l*q^{(\alpha)}*\theta^{(\beta)}(t)\biggr\} da\\
+ & \xi*\int_{B}\biggl\{\chi*\left(\mathcal{L}_{ij}^{(\alpha)}*e_{ij}^{(\beta)}(t)+\mathcal{M}_{ij}^{(\alpha)}*V_{ij}^{(\beta)}(t)\right)\\
+ & \mathcal{R}^{(\alpha)}*\theta^{(\beta)}(t)-\chi*\frac{1}{T_{0}}\mathcal{H}^{(\alpha)}*\theta^{(\beta)}(t)\biggr\} dv,
\end{aligned}
\label{eq:I-cal-def}
\end{equation}
}we have 
\begin{equation}
\mathcal{I}_{\alpha\beta}(t)=\mathcal{I}_{\beta\alpha}(t)\qquad\forall t\in I,\ \forall\alpha,\beta=1,2.\label{eq:I-cal-symmetry}
\end{equation}
\end{thm}

\begin{proof}
We introduce the following function 
\begin{equation}
\mathcal{J}_{\alpha\beta}(t)=\chi*\zeta*\mathcal{T}_{\alpha\beta}(t)-\xi*\chi*\rho\eta^{(\alpha)}*\theta^{(\beta)}(t).\label{eq:J-cal}
\end{equation}
From \eqref{eq:zeta-T} and \eqref{eq:convolution-alternative}$_{8}$
we obtain, with help of eq. \eqref{eq:T-relations},
\begin{equation}
\begin{aligned} & \int_{B}\biggl\{\mathcal{J}_{\alpha\beta}(t)-\chi*\xi*\theta^{0(\alpha)}e^{-t/\beta}*\rho\hat{\eta}^{(\beta)}(t)\\
 & \quad+\xi*\mathcal{R}^{(\alpha)}*\theta^{(\beta)}(t)\biggr\} dv=\\
 & =\mathcal{\int_{B}\biggl\{\chi*\zeta*\hat{T}}_{\alpha\beta}(t)+c\xi*\theta^{(\alpha)}*\theta^{(\beta)}(t)\\
 & \quad-\chi*\xi*\left[\rho\hat{\eta}^{(\alpha)}*\theta^{(\beta)}(t)+\theta^{(\alpha)}*\rho\hat{\eta}^{(\beta)}(t)\right]\biggr\} dv=\\
 & =\int_{B}\biggl\{\mathcal{J}_{\beta\alpha}(t)-\chi*\xi*\theta^{0(\beta)}e^{-t/\beta}*\rho\hat{\eta}^{(\alpha)}(t)\\
 & \quad+\xi*\mathcal{R}^{(\beta)}*\theta^{(\alpha)}(t)\biggr\} dv.
\end{aligned}
\label{eq:J-cal-step-1}
\end{equation}
Using \eqref{eq:J-cal}, \eqref{eq:T-integral}, \eqref{eq:convolution-alternative}$_{3}$,
\eqref{eq:constitutive}$_{6}$ and theorem of divergence, we have{\small{}
\begin{align*}
 & \int_{B}\biggl\{\mathcal{J}_{\alpha\beta}(t)-\chi*\xi*\theta^{0(\alpha)}e^{-t/\beta}*\rho\hat{\eta}^{(\beta)}(t)+\xi*\mathcal{R}^{(\alpha)}*\theta^{(\beta)}(t)\biggr\} dv=\\
 & \ +\chi*\zeta*\int_{B}\biggl[\mathcal{F}_{i}^{(\alpha)}*u_{i}^{(\beta)}(t)-\xi*g^{(\alpha)}*\varphi^{(\beta)}(t)\biggr]dv\\
 & \;+\chi*\xi*\int_{\partial B}\biggl\{\zeta*\biggl[P_{i}^{(\alpha)}*u_{i}^{(\beta)}(t)+R_{i}^{(\alpha)}*\mathcal{D}u_{i}^{(\beta)}(t)\\
 & \;+\Lambda^{(\alpha)}*\varphi^{(\beta)}(t)+H^{(\alpha)}*\mathcal{D}\varphi^{(\beta)}(t)\biggr]+\frac{1}{T_{0}}l*q^{(\alpha)}*\theta^{(\beta)}\biggr\} da\\
 & \ +\xi*\int_{B}\biggl\{\chi*\left(\mathcal{L}_{ij}^{(\alpha)}*e_{ij}^{(\beta)}(t)+\mathcal{M}_{ij}^{(\alpha)}*V_{ij}^{(\beta)}(t)\right)\\
 & \ +\mathcal{R}^{(\alpha)}*\theta^{(\beta)}(t)-\chi*\frac{1}{T_{0}}\mathcal{H}^{(\alpha)}*\theta^{(\beta)}(t)\biggr\} dv\\
 & \ -\chi*\zeta*\int_{B}\rho u_{i}^{(\alpha)}*u_{i}^{(\beta)}(t)dv+\xi*\chi*l*\int_{B}k_{ij}\theta_{,j}^{\mathbf{(\alpha)}}*\theta_{.i}^{(\beta)}(t)dv.
\end{align*}
From this equation and eqs.} \eqref{eq:I-cal-def} 

\begin{align*}
 & \int_{B}\biggl\{\mathcal{J}_{\alpha\beta}(t)-\chi*\xi*\theta^{0(\alpha)}e^{-t/\beta}*\rho\hat{\eta}^{(\beta)}(t)\\
 & \qquad+\xi*\mathcal{R}^{(\alpha)}*\theta^{(\beta)}(t)\biggr\} dv=\\
 & =\mathcal{I}_{\alpha\beta}(t)-\chi*\zeta*\int_{B}\rho u_{i}^{(\alpha)}*u_{i}^{(\beta)}(t)dv\\
 & \qquad+\xi*\chi*l*\int_{B}k_{ij}\theta_{,j}^{\mathbf{(\alpha)}}*\theta_{.i}^{(\beta)}(t)dv,
\end{align*}
so that, with the help of eq. \eqref{eq:J-cal-step-1}, we arrive
to eq. \eqref{eq:I-cal-symmetry}.
\end{proof}

\section{Variational principle}

In this section, we formulate a variational principle for the considered
model. To this aim, we define for each $t\in I$ the functional $\Lambda_{t}$
defined on $\mathcal{V}$ as follows 
\[
\begin{aligned} & \Lambda_{t}\left\{ \pi\right\} =\int_{B}\Biggl\{\chi*\Biggl\{\zeta*\Biggl[\xi*\left(\tau_{ji,j}-\mu_{kji,kj}\right)+\mathcal{F}_{i}\\
 & \;-\frac{1}{2}\rho u_{i}\Biggr]*u_{i}+\zeta*\xi*\left(\sigma_{i,i}-Q_{ji,ji}-g\right)*\varphi\\
 & \:-\xi*\frac{1}{T_{0}}\left(-l*q_{i,i}+\mathcal{H}-\rho T_{0}\eta\right)*\theta\\
 & \:+\xi*\left[\zeta*\left(\tau_{ij}-\frac{1}{2}\hat{\tau}_{ij}\right)+\mathcal{L}_{ij}\right]*e_{ij}\\
 & \:+\zeta*\xi*\left(\mu_{ijk}-\frac{1}{2}\hat{\mu}_{ijk}\right)*\kappa_{ijk}\\
 & \:-\zeta*\xi*\left(\sigma_{i}-\frac{1}{2}\hat{\sigma}_{i}\right)*E_{i}\\
 & \:+\xi*\left[-\zeta*\left(Q_{ij}-\frac{1}{2}\hat{Q}_{ij}\right)+\mathcal{M}_{ij}\right]*V_{ij}\\
 & \:-\xi*\frac{1}{c}\left[-\chi*\frac{1}{2}\left(\rho\eta-\rho\hat{\eta}\right)+\mathcal{R}\right]*\left(\rho\eta-\rho\hat{\eta}\right)\\
 & \:-\xi*l*\left[-\frac{1}{T_{0}}\left(q_{i}-\frac{1}{2}\hat{q}_{i}\right)\right]*\beta_{i}\Biggr\}\Biggr\} dv\\
 & \:-\chi*\zeta*\xi*\left[\int_{S_{1}}P_{i}*\hat{u}_{i}da+\int_{\Sigma_{1}}\left(P_{i}-\hat{P}_{i}\right)*u_{i}da\right]\\
 & \:-\chi*\zeta*\xi*\left[\int_{S_{2}}\!R_{i}*\hat{d}_{i}da\!+\!\int_{\Sigma_{2}}\!\left(R_{i}-\hat{R}_{i}\right)*\mathcal{D}u_{i}da\right]\\
 & \:-\frac{1}{T_{0}}\chi*\xi*l*\left[\int_{S_{3}}q*\hat{\theta}da+\int_{\Sigma_{3}}\left(q-\hat{q}\right)*\theta da\right]\\
 & \:-\chi*\zeta*\xi*\left[\int_{S_{4}}\Lambda*\hat{\varphi}da+\int_{\Sigma_{4}}\left(\Lambda-\hat{\Lambda}\right)*\varphi da\right]\\
 & \:-\chi*\zeta*\xi*\left[\int_{S_{5}}H*\hat{\omega}da+\int_{\Sigma_{5}}\left(H-\hat{H}\right)*\mathcal{D}\varphi da\right].
\end{aligned}
\]
We say that the variation of $\Lambda_{t}$ is zero at $\pi$ over
$\mathcal{V}$ if and only if $\frac{d}{d\lambda}\Lambda_{t}\left\{ \pi+\lambda\pi'\right\} $
exists and is zero for any $\pi'\in\mathcal{V}$ , i.e.
\[
\delta\Lambda_{t}\left\{ \pi\right\} =0\;\Leftrightarrow\;\frac{d}{d\lambda}\Lambda_{t}\left\{ \pi+\lambda\pi'\right\} \Biggr|_{\lambda=0}=0.
\]

\begin{thm}
Fixed $t\in I$, the variation $\delta\Lambda_{t}\left\{ \pi\right\} $
of functional $\Lambda_{t}$ corresponding to $\pi\in\mathcal{V}$
is null if and only if $\pi$ is a solution of the considered mixed
initial-boundary value problem $\Pi$, i.e. $\delta\Lambda_{t}\left\{ \pi\right\} =0$
if and only if $\pi\in\mathcal{K}$. 
\end{thm}

\begin{proof}
To begin, we point out that for any $\pi,\pi'\in\mathcal{V}$ 
\begin{equation}
\begin{aligned} & \int_{\partial B}\left[\left(\tau'_{ji}-\mu'_{kji,k}\right)*u_{i}+\mu'_{jli}*u_{i,l}\right]n_{j}da\\
 & \quad=\int_{\partial B}\left[P'_{i}*u_{i}+R'_{i}*\mathcal{D}u_{i}\right]da,\\
 & \int_{\partial B}\left[\left(\sigma'_{j}-Q'_{ij,i}\right)*\varphi+Q_{ji}*\varphi_{,i}\right]n_{j}da\\
 & \quad=\int_{\partial B}\left[\Lambda'*\varphi+H'*\mathcal{D}\varphi\right]da,
\end{aligned}
\label{eq:variational-step}
\end{equation}
where we take into account eqs.~\eqref{eq:symmetries} and~\eqref{eq:constitutive-reduced}.
Now, by means of the well-known properties of the convolution product,
the definition of variation, eqs.~\eqref{eq:S}, \eqref{eq:S-relations},
\eqref{eq:variational-step} and the divergence theorem, we arrive
to {\small{}
\[
\begin{aligned} & \delta\Lambda_{t}\left\{ \pi\right\} =\int_{B}\chi*\zeta*\left[\xi*\left(\tau_{ji,j}-\mu_{kji,kj}\right)+\mathcal{F}_{i}-\rho u_{i}\right]*u'_{i}\\
 & +\zeta*\xi*\left(\sigma_{i,i}-Q_{ji,ji}-g\right)*\varphi'\\
 & -\xi*\frac{1}{T_{0}}\left(-l*q_{i,i}+\mathcal{H}-\rho T_{0}\eta\right)*\theta'\\
 & +\xi*\zeta*\left[\left(e_{ij}-u_{i,j}\right)*\tau'_{ij}+\left(\kappa_{ijk}-u_{k,ij}\right)*\mu'_{ijk}\right]\\
 & -\xi*\zeta*\left[\left(E_{i}+\varphi_{,i}\right)*\sigma'_{i}+\left(V_{ij}+\varphi_{,ij}\right)*Q'_{ij}\right]\\
 & +\xi*l*\frac{1}{T_{0}}\left[\left(\beta_{i}-\theta_{,i}\right)*q'_{i}+\left(q_{i}-\hat{q}_{i}\right)*\beta'_{i}\right]\\
 & +\xi*\biggl\{\zeta*\left(\tau_{ij}-\hat{\tau}_{ij}\right)+\mathcal{L}_{ij}\\
 & -\frac{1}{c}a_{ij}^{(14)}\biggl[-\chi*\left(\rho\eta-\rho\hat{\eta}\right)+\mathcal{R}\biggr]\biggr\}*e'_{ij}\\
 & +\zeta*\xi*\left[\left(\mu_{ijk}-\hat{\mu}_{ijk}\right)*\kappa'_{ijk}-\left(\sigma_{i}-\hat{\sigma}_{i}\right)*E'_{i}\right]\\
 & +\xi*\biggl\{-\zeta*\left(Q_{ij}-\hat{Q}_{ij}\right)+\mathcal{M}_{ij}\\
 & -\frac{1}{c}a_{ij}^{(47)}\left[-\chi*\left(\rho\eta-\rho\hat{\eta}\right)+\mathcal{R}\right]\biggr\}*V'_{ij}\\
 & -\xi*\frac{1}{c}\left\{ -\chi*\left(\rho\eta-\rho\hat{\eta}\right)+\mathcal{R}^{(4)}-c\theta\right\} *\rho\eta'\Bigr\} dv\\
 & +\chi*\zeta*\xi*\left[\int_{S_{1}}P'_{i}*\left(u_{i}-\hat{u}_{i}\right)da-\int_{\Sigma_{1}}\left(P_{i}-\hat{P}_{i}\right)*u'_{i}da\right]\\
 & +\chi*\zeta*\xi*\left[\int_{S_{2}}R'_{i}*\left(\mathcal{D}u_{i}-\hat{d}_{i}\right)da-\int_{\Sigma_{2}}\left(R_{i}-\hat{R}_{i}\right)*\mathcal{D}u'_{i}da\right]\\
 & +\frac{1}{T_{0}}\chi*\xi*l*\left[\int_{S_{3}}q'*\left(\theta-\hat{\theta}\right)da-\int_{\Sigma_{3}}\left(q-\hat{q}\right)*\theta'da\right]\\
 & \:+\chi*\zeta*\xi*\left[\int_{S_{4}}\Lambda'*\left(\varphi-\hat{\varphi}\right)da-\int_{\Sigma_{4}}\left(\Lambda-\hat{\Lambda}\right)*\varphi'\right]da\\
 & \:+\chi*\zeta*\xi*\left[\int_{S_{5}}H'*\left(\mathcal{D}\varphi-\hat{\omega}\right)da-\int_{\Sigma_{5}}\left(H-\hat{H}\right)*\mathcal{D}\varphi'\right]da,
\end{aligned}
\]
}where $\hat{\tau}'{}_{ij}$, $\hat{\mu}'_{ijk}$, $\hat{\sigma}'_{i}$,
$\hat{Q}'_{ij}$, $\hat{q}'_{i}$, $P'_{i}$, $R'_{i}$, $q'$, $\Lambda'$
and $H'$ are defined by eqs.~\eqref{eq:constitutive-reduced} corresponding
to $\pi'$. Then, for any $\pi'\in\mathcal{V}$ we have that $\delta\Lambda_{t}\left\{ \pi\right\} =0$
if and only if eqs.~\eqref{eq:convolution-alternative}, \eqref{eq:geometric-1},
\eqref{eq:geometric-2}, \eqref{eq:boundary-conditions} hold.
\end{proof}

\section{Conclusions}

In this paper we considered the linear theory of thermopiezoelectric
nonsimple materials as established in \citet{Passarella2022} and
in particular the case of center-symmetric materials. 

We defined a mixed initial-boundary value problem under non-homogeneous
initial conditions and presented a characterization of the mixed initial-boundary
value problem in an alternative way, by including the initial conditions
into the field equations. Starting from a reciprocity relation which
involves two processes at different times, a reciprocity theorem has
been presented. Moreover, a uniqueness result was established without
using the definiteness assumptions on internal energy. Finally, another
reciprocity theorem based on the convolution product and a variational
principle have been derived.

Further developments of this theory could be related to general, non
center-symmetric, materials. Another possible application is to the
study of wave propagation in isotropic materials. 

\bibliographystyle{myunsrtnat}
\bibliography{main}

\begin{thebibliography}{30}
\providecommand{\natexlab}[1]{#1}
\providecommand{\url}[1]{\texttt{#1}}
\expandafter\ifx\csname urlstyle\endcsname\relax
  \providecommand{\doi}[1]{doi: #1}\else
  \providecommand{\doi}{doi: \begingroup \urlstyle{rm}\Url}\fi

\bibitem[Passarella and Tibullo(2022)]{Passarella2022}
F.~Passarella and V.~Tibullo.
\newblock Uniqueness of solutions in thermopiezoelectricity of nonsimple
  materials.
\newblock \emph{submitted to ZAMM}, 2022.

\bibitem[Green and Laws(1972)]{Green1972a}
A.~E. Green and N.~Laws.
\newblock On the entropy production inequality.
\newblock \emph{Arch Ration Mech An}, 45:\penalty0 47--53, 1972.
\newblock ISSN 0003-9527.

\bibitem[Chandrasekharaiah(1998)]{Chandrasekharaiah1998}
D.~S. Chandrasekharaiah.
\newblock Hyperbolic thermoelasticity: {A} review of recent literature.
\newblock \emph{Appl Mech Rev}, 51:\penalty0 705--729, 1998.

\bibitem[Chandrasekharaiah(1986)]{Chandrasekharaiah1986}
D.~S. Chandrasekharaiah.
\newblock Some theorems in generalized micropolar thermoelasticity.
\newblock \emph{Arch Mech}, 38\penalty0 (3):\penalty0 319--328, 1986.
\newblock ISSN 0373-2029.

\bibitem[Ieşan(2004)]{Iesan2004c}
D.~Ieşan.
\newblock \emph{Thermoelastic models of continua}, volume 118 of \emph{Solid
  Mechanics and Its Applications}.
\newblock Springer, Dordrecht, 2004.
\newblock ISBN 978-90-481-6634-3.
\newblock \doi{10.1007/978-1-4020-2310-1}.

\bibitem[Passarella et~al.(2013)Passarella, Tibullo, and
  Zampoli]{Passarella2013}
F.~Passarella, V.~Tibullo, and V.~Zampoli.
\newblock On microstretch thermoviscoelastic composite materials.
\newblock \emph{Eur. J. Mech. A. Solids}, 37:\penalty0 294--303, 2013.
\newblock \doi{10.1016/j.euromechsol.2012.07.002}.

\bibitem[Toupin(1962)]{Toupin1962}
R.~A. Toupin.
\newblock Elastic materials with couple-stresses.
\newblock \emph{Arch Ration Mech An}, 11\penalty0 (1):\penalty0 385--414, 1962.
\newblock \doi{10.1007/BF00253945}.

\bibitem[Toupin(1964)]{Toupin1964}
R.~A. Toupin.
\newblock Theories of elasticity with couple-stress.
\newblock \emph{Arch Ration Mech An}, 17\penalty0 (2):\penalty0 85--112, 1964.
\newblock \doi{10.1007/BF00253050}.

\bibitem[Mindlin(1964)]{Mindlin1964}
R.~D. Mindlin.
\newblock Micro-structure in linear elasticity.
\newblock \emph{Arch Ration Mech An}, 16\penalty0 (1):\penalty0 51--78, 1964.
\newblock \doi{10.1007/BF00248490}.

\bibitem[Toupin and Gazis(1965)]{Toupin1965}
R.~A. Toupin and D.~C. Gazis.
\newblock Surface effects and initial stress in continuum and lattice models of
  elastic crystals.
\newblock In \emph{Lattice Dynamics}, pages 597--605. Elsevier, 1965.

\bibitem[Ahmadi and Firoozbakhsh(1975)]{Ahmadi1975}
G.~Ahmadi and K.~Firoozbakhsh.
\newblock First strain gradient theory of thermoelasticity.
\newblock \emph{Int J Solids Struct}, 11\penalty0 (3):\penalty0 339--345, 1975.
\newblock \doi{10.1016/0020-7683(75)90073-6}.

\bibitem[Batra(1976)]{Batra1976}
R.~C. Batra.
\newblock Thermodynamics of non-simple elastic materials.
\newblock \emph{J Elasticity}, 6\penalty0 (4):\penalty0 451--456, 1976.
\newblock \doi{10.1007/BF00040904}.

\bibitem[Mindlin and Eshel(1968)]{Mindlin1968}
R.~D. Mindlin and N.~N. Eshel.
\newblock On first strain-gradient theories in linear elasticity.
\newblock \emph{Int J Solids Struct}, 4\penalty0 (1):\penalty0 109--124, 1968.
\newblock \doi{10.1016/0020-7683(68)90036-X}.

\bibitem[Ahmadi(1977)]{Ahmadi1977}
G.~Ahmadi.
\newblock Thermoelastic stability of first strain gradient solids.
\newblock \emph{Int. J. Non Linear Mech.}, 12\penalty0 (1):\penalty0 23--32,
  1977.
\newblock \doi{10.1016/0020-7462(77)90013-0}.

\bibitem[Ieşan(1983)]{Iesan1983}
D.~Ieşan.
\newblock Thermoelasticity of nonsimple materials.
\newblock \emph{J Therm Stresses}, 6\penalty0 (2-4):\penalty0 167--188, 1983.
\newblock \doi{10.1080/01495738308942176}.

\bibitem[Ciarletta and Ieşan(1989)]{Ciarletta1989a}
M.~Ciarletta and D.~Ieşan.
\newblock On the nonlinear theory of nonsimple thermoelastic bodies.
\newblock \emph{J Therm Stresses}, 12\penalty0 (4):\penalty0 545--557, 1989.

\bibitem[Kalpakides and Agiasofitou(2002)]{Kalpakides2002}
V.~K. Kalpakides and E.~K. Agiasofitou.
\newblock On material equations in second gradient electroelasticity.
\newblock \emph{J Elasticity}, 67\penalty0 (3):\penalty0 205--227, 2002.
\newblock \doi{10.1023/A:1024926609083}.

\bibitem[Aouadi et~al.(2019)Aouadi, Ciarletta, and Tibullo]{Aouadi2019}
M.~Aouadi, M.~Ciarletta, and V.~Tibullo.
\newblock Analytical aspects in strain gradient theory for chiral {Cosserat}
  thermoelastic materials within three {Green}-{Naghdi} models.
\newblock \emph{J Therm Stresses}, 42\penalty0 (6):\penalty0 681--697, 2019.
\newblock \doi{10.1080/01495739.2019.1571974}.

\bibitem[Aouadi et~al.(2020{\natexlab{a}})Aouadi, Passarella, and
  Tibullo]{Aouadi2020}
M.~Aouadi, F.~Passarella, and V.~Tibullo.
\newblock Exponential stability in {Mindlin's Form II} gradient
  thermoelasticity with microtemperatures of type {III}: {Mindlin's II}
  gradient thermoelastic.
\newblock \emph{Proc. R. Soc. London, Ser. A}, 476\penalty0 (2241),
  2020{\natexlab{a}}.
\newblock \doi{10.1098/rspa.2020.0459}.

\bibitem[Aouadi et~al.(2020{\natexlab{b}})Aouadi, Amendola, and
  Tibullo]{Aouadi2020a}
M.~Aouadi, A.~Amendola, and V.~Tibullo.
\newblock Asymptotic behavior in {Form II Mindlin’s} strain gradient theory
  for porous thermoelastic diffusion materials.
\newblock \emph{J Therm Stresses}, 43\penalty0 (2):\penalty0 191--209,
  2020{\natexlab{b}}.
\newblock \doi{10.1080/01495739.2019.1653802}.

\bibitem[Bartilomo and Passarella(1997)]{Bartilomo1997}
V.~Bartilomo and F.~Passarella.
\newblock Basic theorems for nonsimple thermoelastic solids.
\newblock \emph{Bul. Inst. Politeh. Iaşi. Secţ. I. Mat. Mec. Teor. Fiz.},
  43(47)\penalty0 (1-2):\penalty0 59--70, 1997.

\bibitem[Eringen(2004)]{Eringen2004}
A.~C. Eringen.
\newblock Electromagnetic theory of microstretch elasticity and bone modeling.
\newblock \emph{Int J Eng Sci}, 42\penalty0 (3-4):\penalty0 231--242, 2004.
\newblock \doi{10.1016/S0020-7225(03)00288-X}.

\bibitem[Truesdell and Toupin(1960)]{Truesdell1960}
C.~Truesdell and R.~Toupin.
\newblock The classical field theories.
\newblock In S.~Flügge, editor, \emph{Handbuch der physik}, volume III.
  Springer-Verlag, Berlin - Heidelberg - New York, 1960.

\bibitem[Parkus(1972)]{Parkus1972}
H.~Parkus.
\newblock \emph{Magneto-thermoelasticity}, volume 118.
\newblock Springer, 1972.

\bibitem[Grot(1976)]{Grot1976}
R.~A. Grot.
\newblock Relativistic continuum physics: electromagnetic interactions.
\newblock In A.~C. Eringen, editor, \emph{Continuum physics}, volume III -
  Mixtures and EM Field Theories, pages 129--219. Elsevier, 1976.

\bibitem[Nowacki(1983)]{Nowacki1983}
W.~Nowacki.
\newblock Mathematical models of phenomenological piezo-electricity.
\newblock In R.~H. O.~Brulin, editor, \emph{New Problems in Mechanics of
  Continua}, pages 30--50. University of Waterloo Press, Ontario, 1983.

\bibitem[Maugin(1988)]{Maugin1988}
G.~A. Maugin.
\newblock \emph{Continuum mechanics of electromagnetic solids}, volume~33 of
  \emph{Applied Mathematics and Mechanics}.
\newblock North Holland, Amsterdam, New York, Oxford, Tokio, 1988.

\bibitem[Morro and Straughan(1991)]{Morro1991}
A.~Morro and B.~Straughan.
\newblock A uniqueness theorem in the dynamical theory of piezoelectricity.
\newblock \emph{Math. Methods Appl. Sci.}, 14\penalty0 (5):\penalty0 295--299,
  1991.
\newblock \doi{10.1002/mma.1670140502}.

\bibitem[Passarella et~al.(2011)Passarella, Tibullo, and
  Zampoli]{Passarella2011a}
F.~Passarella, V.~Tibullo, and V.~Zampoli.
\newblock On the heat-flux dependent thermoelasticity for micropolar porous
  media.
\newblock \emph{J Therm Stresses}, 34\penalty0 (8):\penalty0 778--794, 2011.
\newblock \doi{10.1080/01495739.2011.564041}.

\bibitem[Gurtin(1964)]{Gurtin1964}
M.~E. Gurtin.
\newblock Variational principles for linear elastodynamics.
\newblock \emph{Archive for Rational Mechanics and Analysis}, 16\penalty0
  (1):\penalty0 34--50, 1964.
\newblock \doi{10.1007/BF00248489}.

\end{thebibliography}
 
\end{document}